\documentclass[twocolumn]{autart}    

\usepackage{graphicx}          

\usepackage{amssymb}
\usepackage{amsfonts}
\usepackage{mathrsfs}
\usepackage{amsmath}
\usepackage{graphicx}
\usepackage{empheq}
\usepackage{indentfirst}
\usepackage{color}
\usepackage{cases}

\newcommand{\tr}{^{\text{tr}}}

\newcommand{\field}[1]{\ensuremath{\mathbb{#1}}}

\newcommand{\R}{\field{R}\xspace}

\newcommand{\lin}[1]{\mathcal{L}(#1)}

\usepackage{xspace}

\newcommand\Topt{T_{\mathrm{opt}}}
\newcommand\volt{\mathcal{V}}
\newcommand\fred{\mathcal{F}}

\newcommand{\Kone}{K_1}
\newcommand{\Gone}{G_1}
\newcommand{\Gtwo}{G_2}

\newcommand{\ens}[1]{ \left\{#1\right\} }

\newtheorem{theorem}{Theorem}

 \newtheorem{lemma}{Lemma}
 \newtheorem{proposition}{Proposition}
 
 \newtheorem{remark}{Remark}

\begin{document}

\begin{frontmatter}

\title{Finite-time boundary stabilization of  general linear  hyperbolic balance laws via Fredholm backstepping transformation\thanksref{footnoteinfo}}

\thanks[footnoteinfo]{}

\author[Coron]{Jean-Michel Coron}\ead{coron@ann.jussieu.fr},   
\author[Hu]{Long Hu}\ead{hul@sdu.edu.cn},            
\author[Olive]{Guillaume Olive}\ead{math.golive@gmail.com} 

\address[Coron]{Sorbonne Universit\'{e}s, UPMC Univ Paris 06, UMR 7598, Laboratoire
Jacques-Louis Lions, 4 place Jussieu, F-75005 Paris, France.} 
\address[Hu]{School of Mathematics, Shandong University, Jinan, Shandong 250100, China.}             
\address[Olive]{Universit\'{e} de Bordeaux, UMR 5251, Institut de Math\'{e}matiques de Bordeaux, 351 cours de la Lib\'{e}ration, F-33405 Talence, France.}


\begin{abstract}              
This paper is devoted to a simple and new proof on the optimal finite control time for general linear coupled hyperbolic system by using boundary feedback on one side.
The feedback control law is designed by first using a Volterra transformation of the second kind and then using an invertible Fredholm transformation. Both existence and invertibility of the transformations are easily obtained.
\end{abstract}

\end{frontmatter}

\section{Introduction}\label{review linear case}
In this paper, we investigate the stabilization of the following $n \times n$ linear coupled hyperbolic system:
\begin{equation}\label{syst init}
\left\{\begin{array}{l}
u_t(t,x)+\Lambda(x) u_x(t,x)=\Sigma(x) u(t,x), \\
u_-(t,1)=F\left(u(t)\right), \ u_+(t,0)=Q u_-(t,0),  \\
t \in (0,+\infty), \, x\in (0,1),
\end{array}\right.
\end{equation}
where $u=(u_-\tr ,u_+\tr )\tr $ is the state and $F$ is the feedback.
We assume that the matrix $\Lambda \in C^1([0,1])^{n \times n}$ is diagonal: $\Lambda =\mathrm{diag}(\lambda_1,\cdots,\lambda_n)$ and such that $\lambda_i(x) \neq 0$ and $\lambda_i(x) \neq \lambda_j(x)$ for every $x\in [0,1]$, for every $i\in \{1,\ldots,n\}$ and for every 
$j\in \{1,\ldots,n\}\setminus \{i\}$.
Therefore, without loss of generality, we assume that
$$\Lambda=\left(\begin{array}{cc}\Lambda_- & 0 \\0 & \Lambda_+\end{array}\right),$$
where
$$
\Lambda_-=\mathrm{diag}(\lambda_1,\cdots,\lambda_m),
\quad \Lambda_+ = \mathrm{diag} (\lambda_{m+1}, \cdots, \lambda_n),
$$
are diagonal submatrices and
$$\lambda_1(x)<\cdots<\lambda_m(x)<0< \lambda_{m+1}(x)<\cdots<\lambda_n(x),$$
for all $x\in [0,1]$.
Note that we assume that $n\geq 2$ and $m\in\{1,\ldots,n-1\}$.
 Finally, the matrix $\Sigma \in C^0([0,1])^{n \times n}$ couples the equations of the system inside the domain and the constant matrix $Q \in \R^{(n-m) \times m}$ couples the equations of the system on the boundary.

Note that the Riesz representation theorem shows that every bounded linear feedback $F \in \lin{L^2(0,1)^n,\R^m}$ has necessarily the form
\begin{multline}\label{riesz}
Fu=\left(\sum_{j=1}^n \int_0^1 f_{ij}(x) u_j(x) \, dx\right)_{1 \leq i \leq m}, \\ u=(u_1,\ldots,u_n)\tr \in L^2(0,L)^n,
\end{multline}
for some $f_{ij} \in L^2(0,1)$, $i \in \ens{1,\ldots,m}$, $j \in \ens{1,\ldots,n}$.
We can prove that, with this type of boundary conditions, the closed-loop system \eqref{syst init} is well-posed: for every $F \in \lin{L^2(0,1)^n,\R^m}$ and $u^0 \in L^2(0,1)^n$, there exists a unique (weak) solution $u \in C^0([0,+\infty);L^2(0,1)^n)$ to
\begin{equation}\label{Cauchy init}
\left\{\begin{array}{l}
u_t(t,x)+\Lambda(x) u_x(t,x)=\Sigma(x) u(t,x), \\
u_-(t,1)=F\left(u(t)\right), \ u_+(t,0)=Q u_-(t,0),  \\
u(0,x)=u^0(x), \\
t \in (0,+\infty), \, x\in (0,1),
\end{array}\right.
\end{equation}

The purpose of this paper is to find a full-state feedback control law $F$ such that the corresponding closed-loop system \eqref{syst init} vanishes after some time, that is such that there exists $T>0$ such that, for every $u^0 \in L^2(0,1)^n$ for the solution $u \in C^0([0,+\infty);L^2(0,1)^n)$ to \eqref{Cauchy init}, we have
\begin{equation}\label{FTS}
u(t)=0,\quad \forall t \geq T,
\end{equation}
and to obtain the best time $T$ such that \eqref{FTS} holds.

The boundary stabilization problem of 1-D hyperbolic systems have been widely investigated in the literature for almost half a century. The pioneer works date back to \cite{Rauch-Taylor-1974} and \cite{Russell-1978} for linear coupled hyperbolic systems and \cite{Slemord-1983}, \cite{Greenberg-Li-1984}  for the corresponding nonlinear setting, especially for the quasilinear wave equation. For such systems, many articles are based on the  boundary conditions with the following specific form
\begin{align}\label{disspative}
\left(\begin{array}{c}
 u_+(t,0)\\
 u_-(t,1)
 \end{array}\right)
 =\bold{G}\left(\begin{array}{c}
 u_+(t,1)\\
 u_-(t,0)
 \end{array}\right),
\end{align}
where $\bold{G}: \mathbb{R}^{n}\rightarrow \mathbb{R}^n$ is a suitable smooth vector function.
With this boundary condition (\ref{disspative}),  two methods are distinguished to deal with the stability problem of the linear and nonlinear hyperbolic system.
The first one is the so-called characteristic method, which allows us to estimate the related bounds along the characteristic curves.
This method was previously investigated in \cite{Greenberg-Li-1984} for $2\times2$ systems
and in \cite{Qin-1985,Li-1994,Zhao-thesis-1986} for a generalization to $n\times n$ homogeneous nonlinear hyperbolic systems in the framework of $C^1$ norm.
The second one is the control Lyapunov function method, which was introduced in \cite{JMC-Bastin-2014,JMC-Bastin-Andera-2008,Coron-Andrea-Bastin-2007} to analyze the asymptotic behavior of the nonlinear hyperbolic equations in the context of $C^1$ and $H^2$ solutions.
Both of these two approaches guarantee the exponential stability of the nonlinear homogeneous hyperbolic systems provided that the boundary conditions are dissipative to some extent.
Dissipative boundary conditions are standard static boundary output feedback (that is, a feedback of the state values at the boundaries only).
However, there is a drawback of these boundary conditions when inhomogeneous hyperbolic systems are considered, especially the coupling of which are strong enough.
In Section 5.6 of the recent monograph \cite{Bastin-Coron-2016}, the authors provide a counterexample that shows that there exist linear hyperbolic balance laws, which are controllable by open-loop boundary controls, but are impossible to be stabilized under this kind of boundary feedback.

This limitation can be overcome by using the so-called backstepping method, which connects the original system to a target system with desirable stability properties (e.g. exponential stability) via a Volterra transformation the second kind. This method was introduced and developed by M.~Krstic and his co-workers (see, in particular,  the seminal articles \cite{2003-Bokovic-Balogh-Krstic-MCSS,2003-Liu-SICON,2004-Smyshlyaev-Krstic-IEEE} and the tutorial book \cite{KSbook}).
In \cite{Coron-Vazquez-Krstic-Bastin-SICON}, the authors designed a full-state feedback
control law, with actuation on only one side of the boundary, in order to achieve $H^2$ exponential stability of the closed-loop $2\times 2$ quasilinear hyperbolic system by using Volterra-type backstepping transformation.
Moreover, with this method we can even steer the corresponding linearized hyperbolic system to rest in finite time, that is what is called finite time stabilization.
The presented method can also be extended to linear systems with only one negative characteristic velocity (see \cite{Meglio}).
In \cite{Hu-Dimeglio-Vazquez-Krstic-2015-linear}, a fully general case of coupled heterodirectional hyperbolic PDEs, allowing an arbitrary number of PDEs convecting in each direction and the boundary controls applied on one side, is presented. The proposed boundary controls also yield the finite-time convergence to zero with the control time given by
 \begin{align}\label{t_F}
t_F=\int_0^1 \frac{1}{\lambda_{m+1}(x)} \, dx+\sum_{i=1}^m \int_0^1 \frac{1}{|\lambda_i(x)|} \, dx.
\end{align}
However, this time $t_F$ is larger than the theoretical optimal one we expect and that is given in \cite{Li-Rao-CAM-2010}, namely
\begin{align}\label{vanishing time}
\Topt=\int_0^1 \frac{1}{\lambda_{m+1}(x)} \, dx+\int_0^1 \frac{1}{|\lambda_{m}(x)|} \, dx.
\end{align}
In \cite{Auriol-Dimeglio-2015}, the authors found a minimum time stabilizing controller which makes the coupled hyperbolic system (\ref{syst init}) with constant coefficients vanishes after $\Topt$ by slightly changing the target system in \cite{Hu-Dimeglio-Vazquez-Krstic-2015-linear}, in which only local cascade coupling terms are involved in the PDEs.

In this paper, we show that this kind of controller can be established in a much easier way.
Inspired by the known results of \cite{Hu-Dimeglio-Vazquez-Krstic-2015-linear} and \cite{Hu-Vazquez-Meglio-Krstic-2015}, we will map the initial coupled hyperbolic system (\ref{syst init}) to a new target system in which the cascade coupling terms of the previous works (namely, $G(x)\beta(t,0)$ in \cite{Hu-Dimeglio-Vazquez-Krstic-2015-linear} and $\Omega(x)\beta(t,x)$ in \cite{Auriol-Dimeglio-2015}) can be completely cancelled.
Our strategy is to first transform (\ref{syst init}) to the target system of \cite{Hu-Vazquez-Meglio-Krstic-2015} by a Volterra transformation of the second kind, which is always invertible if the kernel belongs to $L^2$.
Then, regarding the target system obtained as the initial hyperbolic system to be studied, by using a Fredholm transformation as introduced in \cite{Coron-Lu-2014}, we then map this intermediate system to a new target system, vanishing after $\Topt$, without any coupling terms in the PDEs other than a simple trace coupling term.
Moreover, the existence and the invertibility of such a transformation will be easily proved (we point out here that these transformations are not always invertible, see \cite{Coron-Hu-Olive-2015}, but this will indeed be the case here thanks to the cascade structure of the kernel involved in our Fredholm transformation).
Finally, the target system and the original system share the same stability properties due to the invertibility of the transformation.

The main result of this paper is the following:
\begin{theorem}\label{main thm}
There exists $F \in \lin{L^2(0,1)^n,\R^m}$ such that, for every $u^0 \in L^2(0,1)^n$, the solution $u \in C^0([0,+\infty);L^2(0,1)^n)$ to \eqref{Cauchy init} satisfies
$$u(t)=0,\quad \forall t \geq \Topt,$$
where $\Topt$ is given by \eqref{vanishing time}.
\end{theorem}
\begin{remark}
We recall that this result has already been obtained in \cite{Auriol-Dimeglio-2015} in the case of constant matrices $\Lambda$ and $\Sigma$.
Therefore, Theorem \ref{main thm} generalizes this result.
We also believe that, even in the case of constant matrices, the approach we shall present below, based on an invertible Fredholm transformation and a simple target system, is easier than the one presented in \cite{Auriol-Dimeglio-2015}, where a Volterra transformation and a different target system are used. In particular,   we do not need repeatedly use the successive approximation approach to find the kernels in the transformation.
\end{remark}

The rest of the paper is organized as follows.
In Section \ref{sec:targetsystem}, we first recall the results of \cite{Hu-Vazquez-Meglio-Krstic-2015} and then we present a new target system which vanishes after the optimal time $\Topt$.
Then, in Section \ref{sec:backstepping}, we prove the existence of an invertible Fredholm transformation that maps the target system introduced in \cite{Hu-Vazquez-Meglio-Krstic-2015} into the new designed target system.

\section{New target system}\label{sec:targetsystem}

In \cite[Section 2.1]{Hu-Vazquez-Meglio-Krstic-2015} the authors introduced the following target system in the particular case $H=0$:
\begin{equation}\label{w syst}
\left\{\begin{array}{l}
\gamma_t(t,x)+\Lambda(x) \gamma_x(t,x)=G(x) \gamma(t,0), \\
\gamma_-(t,1)=H(\gamma(t)), \ \gamma_+(t,0)=Q \gamma_-(t,0), \\
t \in(0,+\infty), \, x \in (0,1),
\end{array}\right.
\end{equation}
where $\gamma=(\gamma_-\tr ,\gamma_+\tr )\tr $ is the state and $H$ is a feedback.
The matrix $G \in L^{\infty}(0,1)^{n \times n}$ is a lower triangular matrix with the following structure
\begin{align}\label{matrixG}
G=\left(\begin{array}{cc}\Gone & 0 \\\Gtwo & 0\end{array}\right),
\end{align}
where $\Gone \in L^{\infty}(0,1)^{m \times m}$ has the cascade structure
\begin{align}\label{struct g1}
	\Gone&=\begin{pmatrix}
		0&\cdots&\cdots&0\\
		g_{2 \, 1}&\ddots&\ddots&\vdots\\
		\vdots&\ddots&\ddots&\vdots\\
		g_{m \, 1}&\cdots&g_{m \, m-1}&0
	\end{pmatrix},
\end{align}
for some $g_{ij} \in L^{\infty}(0,1)$, $i \in \ens{2,\ldots,m}$, $\, j\in \ens{1, \ldots, i-1}$, and $\Gtwo \in L^{\infty}(0,1)^{(n-m) \times m}$.
We recall that, for every $H \in \lin{L^2(0,1)^n,\R^m}$ and $\gamma^0 \in L^2(0,1)^n$, there exists a unique (weak) solution $\gamma \in C^0([0,+\infty);L^2(0,1)^n)$ to \eqref{w syst} satisfying $\gamma(0,\cdot)=\gamma^0$.

Taking into account the form of the feedbacks (see \eqref{riesz}) we can use the standard backstepping method and establish the following result, in the exact same way as it was done in \cite{Hu-Vazquez-Meglio-Krstic-2015} for the case $H=0$:
\begin{lemma}\label{main lem}
There exist $G \in L^{\infty}(0,1)^{n \times n}$ with the structure \eqref{matrixG}-\eqref{struct g1} and an invertible bounded linear map $\volt:L^2(0,1)^n \longrightarrow L^2(0,1)^n$ such that, for every $H \in \lin{L^2(0,1)^n,\R^m}$, there exists $F \in \lin{L^2(0,1)^n,\R^m}$ such that, for every $u^0 \in L^2(0,1)^n$, if $\gamma \in C^0([0,+\infty),L^2(0,1)^n)$ denotes the solution to
\eqref{w syst} satisfying the initial data $\gamma(0,\cdot)=\volt^{-1}u^0$, then
$$u(t)=\volt \gamma(t),$$
is the solution to the Cauchy problem \eqref{Cauchy init}.
\end{lemma}

For the rest of the paper, $G$ is fixed as in Lemma \ref{main lem}.

In \cite{Hu-Vazquez-Meglio-Krstic-2015}, the authors chose the simplest possibility $H=0$ so that, due to the cascade structure \eqref{matrixG}-\eqref{struct g1}, any solution to the resulting system \eqref{w syst} defined at time $0$ vanishes after the time $t_{F}$ given by (\ref{t_F}) (see \cite[Proposition 2.1]{Hu-Vazquez-Meglio-Krstic-2015} for more details).
However, this appears to be not the best choice since it does not give the expected optimal time $\Topt$.
In the present paper, we will show how to properly choose $H$ in order to reduce the vanishing time to $\Topt$.
For this purpose, the idea is to apply a second time the backstepping method and find a Fredholm mapping that transforms the previous target system \eqref{w syst} into the following new target system:
\begin{equation}\label{g syst}
\left\{\begin{array}{l}
z_t(t,x)+\Lambda(x) z_x(t,x)=\widetilde{G}(x) z(t,0), \\
z_-(t,1)=0, \quad z_+(t,0)=Q z_-(t,0),
\\
t \in(0,+\infty), \, x \in (0,1),
\end{array}\right.
\end{equation}
where $z=(z_-\tr,z_+\tr)\tr$ is the state and $\widetilde{G} \in L^{\infty}(0,1)^{n \times n}$ is the following matrix
\begin{align}\label{Matrix:M}
\widetilde{G}(x)=\left(\begin{array}{cc}0 & 0 \\\Gtwo(x) & 0\end{array}\right),
\end{align}
where $\Gtwo$ is defined in (\ref{matrixG}).
We recall that, for every $z^0 \in L^2(0,1)^n$, there exists a unique (weak) solution $z \in C^0([0,+\infty);L^2(0,1)^n)$ to \eqref{g syst} satisfying $z(0,\cdot)=z^0$.
Moreover one has the following proposition:

\begin{proposition}\label{pro2.1}
For every $z^0 \in L^2(0,1)^n$, the solution $z \in C^0([0,+\infty);L^2(0,1)^n)$ to \eqref{g syst} satisfying $z(0,\cdot)=z^0$ verifies $z(t)=0$ for every $t \geq \Topt$.
\end{proposition}

\emph{Proof.}
Indeed,  using the method of characteristics and the cascade structure \eqref{Matrix:M} of $\widetilde{G}$, one first gets that $z_-(t)=0$ for $t\geq \int_0^1 1/|\lambda_{m}(x)|\, dx$ and then that $z_+(t)=0$ for $t\geq \Topt$.  \qed

We will prove the following result:

\begin{proposition}\label{main prop}
There exist an invertible bounded linear map $\fred:L^2(0,1)^n \longrightarrow L^2(0,1)^n$ and $H \in \lin{L^2(0,1)^n,\R^m}$ such that, for every $\gamma^0 \in L^2(0,1)^n$, if $z \in C^0([0,+\infty),L^2(0,1)^n)$ denotes the solution to  \eqref{g syst} satisfying the initial data $z(0,\cdot)=\fred^{-1}\gamma^0$, then
$$\gamma(t)=\fred z(t),$$
is the solution to \eqref{w syst} satisfying $\gamma(0,\cdot)=\gamma^0$.
\end{proposition}

\begin{remark}
In Lemma \ref{main lem} it is showed that we can reach system \eqref{syst init} from system \eqref{w syst} whatever the feedback $H$ is, $F$ being fixed consequently.
Note that there is no such freedom in Proposition \ref{main prop} as we need the boundary condition $z_-(t,1)=0$ in a crucial way for the proof, see \eqref{cond zu} below.
\end{remark}

Combining all the aforementioned results, it is now easy to obtain Theorem \ref{main thm}:

\emph{Proof of Theorem \ref{main thm}.}
Let $\fred$ and $H$ be the two mappings provided by Propositon \ref{main prop} and then let $\volt$ and $F$ be the corresponding mappings provided by Lemma \ref{main lem}.
Let $z \in C^0([0,+\infty),L^2(0,1)^n)$ be the solution to  \eqref{g syst} associated with the initial data $z(0,\cdot)=\left(\volt \circ \fred\right)^{-1}u^0$.
Then,
\begin{equation}\label{u g}
u(t)=\volt \circ \fred z(t),
\end{equation}
is the solution to the Cauchy problem \eqref{Cauchy init}.
By Proposition \ref{pro2.1}, we know that $z(t)=0$ for every $t \geq \Topt$ and it readily follows from \eqref{u g} that $u(t)=0$ for every $t \geq \Topt$ as well.\qed

Therefore, it only remains to establish Proposition \ref{main prop}.
This is achieved in the next section.

\section{Existence of an invertible Fredholm transformation}\label{sec:backstepping}

In this section we prove Proposition \ref{main prop}.
To this end, we look for a Fredholm transformation $\fred:L^2(0,1)^n \longrightarrow L^2(0,1)^n$:
\begin{multline}\label{def f transfo}
\fred z(x)=z(x)-\int_{0}^1 K(x,y) z(y)dy, \\ x \in (0,1), \, z \in L^2(0,1)^n,
\end{multline}
with a kernel $K \in L^2((0,1)\times(0,1))^{n \times n}$ with the following structure:
\begin{align}\label{Matrix:K}
K=\left(\begin{array}{cc} \Kone & 0 \\ 0 & 0\end{array}\right),
\end{align}
in which $\Kone\in L^2((0,1)\times(0,1))^{m \times m}$ is a lower triangular matrix with 0 diagonal entries, that is has the following cascade structure
\begin{align}\label{2.12s}
	\Kone&=\begin{pmatrix}
		0&\cdots&\cdots&0\\
		k_{2 \, 1}&\ddots&\ddots&\vdots\\
		\vdots&\ddots&\ddots&\vdots\\
		k_{m \, 1}&\cdots&k_{m \, m-1}&0
	\end{pmatrix},
\end{align}
for some $k_{ij} \in L^2((0,1)\times(0,1))$, $i \in \ens{2,\ldots,m}$, $j\in \ens{1,\ldots,i-1}$, yet to be determined.
Note that $\fred$ is clearly invertible due to this very particular structure (see the Appendix \ref{invertibility} for details).
Therefore, we only have to check that $\gamma$ defined by
\begin{equation}\label{w fred}
\gamma(t,x)=z(t,x)-\int_0^1 K(x,y) z(t,y) \, dy,
\end{equation}
is solution to \eqref{w syst} for some $H \in \lin{L^2(0,1)^n,\R^m}$ to be determined as well.

Let us first perform some formal computations to derive the equations that the $k_{ij}$ have to satisfy.
Taking the derivative with respect to time in \eqref{w fred}, using the equation satisfied by $z$ (see the first line of \eqref{g syst}) and integrating by parts yield
\begin{equation*}
\begin{aligned}
& \gamma_t(t,x) =  z_t(t,x)-\int_0^1 K(x,y) z_t(t,y) \, dy\\
 &=-\Lambda(x) z_x(t,x)+\widetilde{G}(x)z(t,0) \\
&\quad +\int_0^1 K(x,y) \Lambda(y) z_y(t,y) \, dy\\
&\quad - \int_0^1 K(x,y)\widetilde{G}(y)z(t,0) \, dy\\
 &=-\Lambda(x) z_x(t,x)+\widetilde{G}(x)z(t,0)
+K(x,1)\Lambda(1) z(t,1)\\
&\quad - K(x,0)\Lambda(0) z(t,0)  -\int_0^1 K_y(x,y) \Lambda(y) z(t,y) \, dy \\
&\quad -\int_0^1 K (x,y)\Lambda_y(y)z(t,y) \, dy \\
& \quad - \int_0^1 K(x,y)\widetilde{G}(y)z(t,0) \, dy.
\end{aligned}
\end{equation*}

Now observe that, since $z_-(t,1)=0$ and because of the structures of $K$ (see (\ref{Matrix:K})) and $\widetilde{G}$ (see (\ref{Matrix:M})), we have the following two conditions:
\begin{equation}\label{cond zu}
K(x,1)\Lambda(1)z(t,1)=0,
\end{equation}
$$K(x,y)\widetilde{G}(y)=0.$$
Therefore,
\begin{align*}
\begin{aligned}
& \gamma_t(t,x)=
-\Lambda(x) z_x(t,x)+\Big(\widetilde{G}(x)- K(x,0)\Lambda(0)\Big) z(t,0)  \\
& \quad -\int_0^1 \Big(K_y(x,y) \Lambda(y)+K (x,y)\Lambda_y(y) \Big)z(t,y) \, dy.
\end{aligned}
\end{align*}
On the other hand, taking the derivative with respect to space in \eqref{w fred} we have
  \begin{align*}
 \begin{aligned}
\gamma_x(t,x)&=z_x(t,x)-\int_0^1 K_x(x,y) z(t,y) \, dy.
 \end{aligned}
 \end{align*}
As a result, we obtain
\begin{align*}
\begin{aligned}
& \gamma_t(t,x)+ \Lambda(x) \gamma_x(t,x)-G(x)\gamma(t,0) \\
& =\Big(\widetilde{G}(x)- K(x,0)\Lambda(0)-G(x)\Big) z (t,0) \\
& \quad -\int_0^1 \Big(K_y(x,y) \Lambda(y)+K (x,y)\Lambda_y(y)\\
& \quad +\Lambda(x) K_x(x,y)-G(x)K(0,y))\Big) z(t,y) \, dy,
\end{aligned}
\end{align*}
and the right-hand side has to be zero.
This yields to the following kernel system
\begin{multline*}
\Lambda(x) K_x(x,y)+K_y(x,y) \Lambda(y)\\+K (x,y)\Lambda_y(y)
-G(x)K(0,y)=0
\end{multline*}
with the condition
$$K(x,0)=(\widetilde{G}(x)-G(x))\Lambda^{-1}(0).$$
In order to guarantee the well-posedness of the system satisfied by $K$, we impose the following extra condition:
\begin{equation}\label{art cond}
K_-(0,y)=0,
\end{equation}
(where $K_-$ denotes the submatrix containing the first $m$ rows of $K$), which turns out to also imply the following, because of the structures of $G$ (see \eqref{matrixG}) and $K$ (see \eqref{Matrix:K}),
$$G(x)K(0,y)=0,$$
and therefore makes the kernel system much simpler to solve.
To summarize, $K$ will satisfy the system
\begin{equation*}
\left\{\begin{array}{l}
\Lambda(x) K_x(x,y)+K_y(x,y) \Lambda(y)+K (x,y)\Lambda_y(y)=0, \\
K_-(0,y)=0, \\
K(x,0)=(\widetilde{G}(x)-G(x))\Lambda^{-1}(0), \\
x,y \in (0,1).
\end{array}\right.
\end{equation*}
Note that the structure \eqref{Matrix:K} of $K$, \eqref{w fred} and \eqref{art cond} imply that
$$\gamma(t,0)=z(t,0).$$
Therefore, the boundary condition at $x=0$ for $\gamma$ is automatically guaranteed:
$$\gamma_+(t,0)=z_+(t,0)=Qz_-(t,0)=Q\gamma_-(t,0).$$
Now, because of the structures of $K$, $\widetilde{G}$ and $G$ given in 
\eqref{Matrix:K}, \eqref{Matrix:M} and \eqref{matrixG} respectively, the system for $K$ translates into the following system for $\Kone$:
\begin{equation*}
\left\{\begin{array}{l}
\Lambda_-(x) (\Kone)_x(x,y)+(\Kone)_y(x,y) \Lambda_-(y)
\\
\phantom{\Lambda_-(x) (\Kone)_x(x,y)+bbb}+\Kone(x,y)(\Lambda_-)_y(y)=0, \\
\Kone(0,y)=0, \\
\Kone(x,0)=-\Gone(x)\Lambda_-^{-1}(0), \\
x,y \in (0,1).
\end{array}\right.
\end{equation*}
Regarding $y$ as the time parameter, this is a standard time-dependent uncoupled hyperbolic system with only positive speeds $\lambda_i(x)/\lambda_j(y)>0$, $i,j \in \ens{1,\ldots,m}$, and therefore it admits a unique (weak) solution $\Kone \in L^2((0,1)\times(0,1))^{m \times m}$.
Actually, using the method of characteristics, we see that the solution is explicitely given by
\begin{equation}\label{explicit}
k_{ij}(x,y)=\frac{g_{ij}\big(\phi_i^{-1}\big(\phi_i(x)-\phi_j(y)\big)\big)}{-\lambda_j(y)},
\end{equation}
if $i \in \ens{2,\ldots,m}$, $j \in \ens{1,\ldots,i-1}$ and $\phi_i(x) \leq \phi_j(y)$, and $k_{ij}(x,y)=0$ otherwise, where
$$\phi_i(x)=\int_0^x \frac{1}{\lambda_i(\xi)}d\xi,\quad i \in \ens{1,\ldots,m}.$$
Note that $\phi_i$ is indeed invertible since it is a monotonically decreasing continuous function of $x$.
Finally, we readily see from \eqref{explicit} that
$$\Kone(1,\cdot) \in L^2(0,1)^{m \times m},$$
so that the map $H:L^2(0,1)^n \longrightarrow \R^m$ given by
$$H\gamma=-\int_0^1 \Kone(1,y)\left[\fred^{-1}\gamma\right]_-(y) \, dy, \quad \gamma \in L^2(0,1)^n,$$
is well-defined and $H \in \lin{L^2(0,1)^n,\R^m}$.
This concludes the proof of Proposition \ref{main prop}.
\qed

\begin{remark}
Let us conclude this paper by pointing out that it would be very interesting to know the target systems that can be achieved with general linear transformations.
We recall that it is proved in \cite{2015-Coron-ICIAM} that, for the finite dimensional control system $\dot y =Ay +Bu$, the target system $\dot y =Ay-\lambda y+Bu$ can be achieved by a linear transformation for every $\lambda \in \mathbb{R}$, if we assume that it is controllable (which is a necessary condition to the rapid stabilization). 
\end{remark}

\begin{ack}                       
The authors thank Amaury Hayat and Shengquan Xiang for useful comments.
This project was  supported by the ERC advanced grant 266907 (CPDENL) of the 7th Research Framework Programme (FP7), ANR Project Finite4SoS (ANR 15-CE23-0007), the Young Scholars
Program of Shandong University (No. 2016WLJH52), the Natural Science Foundation of China (No. 11601284) and the China Postdoctoral Science Foundation (No. BX201600096).   
\end{ack}

\bibliographystyle{plain}       
\bibliography{autosam}

\appendix
\section{Invertibility of the Fredholm transformation}\label{invertibility}

For the completeness we prove in this appendix the invertibility of the Fredholm transformation $\fred$.

\begin{lemma}\label{appendix lemma}
For any given $K\in L^2((0,1)\times(0,1))^{n \times n}$ with the cascade structure (\ref{Matrix:K})-(\ref{2.12s}), the transformation $\fred$ defined by \eqref{def f transfo} is invertible.
Moreover, its inverse has the same form:
\begin{multline*}
\fred^{-1} \gamma(x)=\gamma(x)-\int_{0}^1 \widetilde{\Theta}(x,y) \gamma(y)dy, \\ x \in (0,1), \, \gamma \in L^2(0,1)^n,
\end{multline*}
for some $\widetilde{\Theta}\in L^2((0,1)\times(0,1))^{n\times n}$ with the same structure as $K$, that is,
\begin{align*}
\widetilde{\Theta}=\left(\begin{array}{cc}\Theta & 0 \\ 0 & 0\end{array}\right),
\end{align*}
in which $\Theta\in L^2((0,1)\times(0,1))^{m \times m}$ is a lower triangular matrix with 0 diagonal entries as $\Kone$:
\begin{align*}
	\Theta&=\begin{pmatrix}
		0&\cdots&\cdots&0\\
		\theta_{2 \, 1}&\ddots&\ddots&\vdots\\
		\vdots&\ddots&\ddots&\vdots\\
		\theta_{m \, 1}&\cdots&\theta_{m \, m-1}&0
	\end{pmatrix},
\end{align*}
for some $\theta_{ij} \in L^{2}((0,1)\times(0,1))$, $i \in \ens{2,\ldots,m}$ and $j\in \ens{1, \ldots, i-1}$.
\end{lemma}
\emph{Proof of Lemma \ref{appendix lemma}.}
Let $\gamma=\fred z$, where $z \in L^2(0,1)^n$ is given.
Thanks to \eqref{Matrix:K} and \eqref{w fred}, we have
$$z_i=\gamma_i, \quad \forall i \in \ens{m+1,\ldots,n}.$$
On the other hand, thanks to \eqref{2.12s} and \eqref{w fred}, we have
$$
\left\{\begin{array}{l}
\gamma_1=z_1, \\
\gamma_i=z_i-\sum_{j=1}^{i-1} \int_0^1 k_{ij}(\cdot,y) z_j(y) \, dy, \, \forall i \in \ens{2,\ldots,m}.
\end{array}
\right.
$$
By induction we readily see that
$$
\left\{\begin{array}{l}
z_1=\gamma_1, \\
z_i=\gamma_i-\sum_{j=1}^{i-1} \int_0^1 \theta_{ij}(\cdot,y) \gamma_j(y) \, dy, \, \forall i \in \ens{2,\ldots,m},
\end{array}
\right.
$$
for some $\theta_{ij} \in L^{2}(0,1)$ depending only on $k_{p j}$ for $p \in \ens{j+1,\ldots,i}$.
This proves Lemma \ref{appendix lemma}.
\qed

\end{document}